\documentclass{article}

\usepackage{amsmath,amssymb,latexsym,epsfig}

\usepackage{amsfonts}

\newtheorem{theorem}{Theorem}[section]
\newtheorem{proposition}[theorem]{Proposition}
\newtheorem{corollary}[theorem]{Corollary}
\newtheorem{definition}[theorem]{Definition}
\newtheorem{lemma}[theorem]{Lemma}

\newcommand{\beha}{\begin{enumerate}}
\newcommand{\behe}{\end{enumerate}}

\newcommand{\bR}{{\bf R}}
\newcommand{\bC}{{\bf C}}
\newcommand{\bZ}{{\bf Z}}
\newcommand{\bN}{{\bf N}}

\newcommand{\RE}[1]{(\ref{#1})}
\newcommand{\eproof}{$\quad\Box$\\[0.3cm]}

\def\cal#1{\fam2#1}

\begin{document}
\title{Hausdorff and topological dimension for polynomial automorphisms of $\bC^2$}

\renewcommand{\baselinestretch}{1.0}
\author{
Christian Wolf \footnote{The author was supported by a research
fellowship of the Deutsche Forschungsgemeinschaft (DFG)}
\\
 Department of Mathematics, Indiana University\\
 Rawles Hall, Bloomington, IN 47405, USA\\
 Email: cwolf@indiana.edu}
 \maketitle
\begin{abstract}Let $g$ be a polynomial automorphism
of $\bC^2$. We  study the Hausdorff dimension and  topological
dimension of the Julia set of $g$. We show that when $g$ is a
hyperbolic mapping, then   the Hausdorff dimension of the Julia
set is strictly greater than its topological dimension. Moreover,
the Julia set  cannot be locally connected. We also provide
estimates for the dimension of the Julia sets   in the general
(not necessarily hyperbolic) case.
\end{abstract}
\section{Introduction}
Let $g$ be a polynomial automorphism of $\bC^2$. We can associate
with $g$ a dynamical  degree $d$ which is a conjugacy invariant,
see Section $2$ for the definition. We are interested in
nontrivial dynamics which can occur only for $d>1$. We define
$K^\pm$ as the set of points in $\bC^2$ with bounded
forward/backward orbits, $K=K^+\cap K^-$, $J^\pm=
\partial K^\pm$ and $J= J^+\cap J^-$.
 We refer to $J^\pm$ as the
positive/negative Julia set of $g$ and $J$ is the Julia set of
$g$. The set $J^\pm$ is unbounded and connected, while $K$  and
$J$ are compact. The "chaotic" dynamic (recurrent dynamics with
sensitive dependence on the initial conditions in forward and in
backward time) can only occur on the Julia set $J$. The
nontriviality  of the dynamics of $g$ is reflected by the fact
that $g|_K$ has positive topological entropy equal to $\log d$.
Note that the complex Jacobian determinant $\det\, Dg$ is constant
in $\bC^2$. We will use $a=\det\, Dg$ as standard notation in this
paper. We can restrict our considerations to the volume decreasing
case $(|a|<1)$ and to the volume preserving case $(|a|=1)$,
because otherwise we can consider $g^{-1}$. In this paper we will
use $|a|\leq 1$ as a standing assumption.\\[0.3cm] An important
and nontrivial feature of a dynamical system is the Hausdorff
dimension of its closed invariant sets. For a polynomial map of
$\bC$ this subject has been studied successfully, especially in
the case of hyperbolic and parabolic maps, see \cite{U} for a
recent survey article. The central part of the theory is based on
the so-called Bowen- Ruelle formula, which determines the
Hausdorff dimension of the Julia set as the zero of the pressure
function.\\[0.3cm] In the case of polynomial automorphisms of
$\bC^2$ the situation is more complicated because in contrast to
polynomials of $\bC$, the derivative $Dg(p)$ induces  a
contracting and an expanding direction in the tangent space.
Nevertheless, when $g$ is a hyperbolic mapping,  we also have a
version of the Bowen- Ruelle formula; it provides the Hausdorff
dimension of the unstable/stable slice in terms of the zero of the
unstable/stable pressure function. We will use this result to
analyze the Hausdorff dimension of $J$.
\\[0.3cm] Let us now describe our results more precisely. \\[0.3cm]
First,  we consider the topological dimension of the Julia set of
a hyperbolic mapping $g$. In particular, we show that the
topological dimension of the Julia set can be expressed in terms
of an intersection formula, that is,
\begin{equation}
dim_{top} J = dim_{top} J^+ + dim_{top} J^- - 4.
\end{equation}
The proof of this result is based on general features about
hyperbolic sets. Note that the analogous intersection formula for
the Hausdorff dimension also holds (see \cite{W3}). As an
immediate consequence the topological dimension of $J$ can be at
most $2$. We also provide a criterion which implies that $J$ has
topological dimension at least $1$. Finally, we  show that $J$
cannot be locally connected.\\[0.3cm] Second, we deal with the
Hausdorff dimension of the Julia set of a hyperbolic mapping $g$.
Let $P^{u/s}$ denote the unstable/stable pressure function of $g$
(see Section $4$ for the definition). Then the result of Verjovsky
and Wu \cite{VW} implies that $t^{u/s}=dim_H W^{u/s}_\epsilon (p)
\cap J$ is given by the unique solution of the Bowen- Ruelle
formula
\begin{equation}
P^{u/s}(t)=0.
\end{equation}
Moreover,  $dim_H J=t^u+t^s$ and $dim_H J^\pm = t^{u/s}+2$, see
\cite{W3} and the references therein. We show that the pressure
functions are related to each other in terms of the Jacobian
determinant. As a consequence we obtain the  inequality
\begin{equation}
t^s\leq \frac{t^u \log d}{\log d -t^u \log |a|}.
\end{equation} In particular, $t^s< t^u$ for $|a|<1$ and $t^u=t^s$
for $|a|=1$. This implies that when $|a|\leq d^{-1/2}$, then the
stable slice $W^s_\epsilon(p)\cap J$ is totally disconnected.
Therefore, $dim_{top} J=2$ can occur only for $|a|$ being
sufficiently close to $1$. Combining these results with a
perturbation argument yields.
\\[0.3cm]{ Theorem 4.11.} {\it If $g$ is hyperbolic, then}
$dim_{top}J<dim_H J$.\\[0.3cm] This property is sometimes
(especially in applied sciences) used for  defining  fractal
sets.\\[0.3cm] Last, we consider the dimension of the Julia sets
of a general (not necessarily hyperbolic) polynomial automorphism
$g$. We show that when $g$ is not volume preserving, then the
upper box-dimension of $K$ is strictly smaller than $4$. We also
obtain that if $|a|$ is small, then the Hausdorff dimension of
$J^-$ is close to 2. This generalizes a result of Forn{\ae}ss and
Sibony \cite{FS} about complex H\'enon mappings which are small
perturbations of hyperbolic quadratic  polynomials in $\bC$.
\\[0.3cm]
 This paper
is organized as follows. In Section 2 we present the basic
definitions and notations. The results about the topological
dimension of the Julia set are presented in Section 3. Section  4
is devoted to the analysis of the Hausdorff dimension of the Julia
set. In particular, we show that the Hausdorff dimension of the
Julia set always exceeds its topological dimension. In Section 5
we derive Hausdorff dimension estimates for general polynomial
automorphisms of $\bC^2$.
\section{Notation and Preliminaries}
Let $g$ be a polynomial automorphism of $\bC^2$. We denote by $deg
(g)$ the maximum of the algebraic degree of the components of $g$.
The dynamical degree of $g$ is defined by
\begin{displaymath}
d=\lim\limits_{n \to\infty} \left(deg (g^n)\right)^{1/n}
\end{displaymath}
(see \cite{BS2}, \cite{FM}). According to the Friedland-Milnor
classification Theorem \cite{FM}, each polynomial automorphism of
$\bC^2$ is either conjugate to an elementary polynomial
automorphism of $\bC^2$ (with dynamical degree 1 and trivial
dynamics) or to a finite composition of generalized H\'enon
mappings (with dynamical degree greater than 1 and nontrivial
dynamics).   A finite composition of generalized H\'enon mappings
is a mapping
\begin{equation}\label{defh1}
g=g_1\circ...\circ g_m, \end{equation} where each $g_i$ has the
form
\begin{equation}\label{defhen}
g_i(z,w)= (w,P_i(w)+a_i z),
\end{equation}
$P_i$ is a complex polynomial of degree  $d_i\geq 2$ and $a_i$ is
a non-zero complex number. Note that the dynamical degree of $g$
is equal to $d=d_1\cdot...\cdot d_m$ and therefore coincides with
$deg(g)$. The complex Jacobian determinant of $g$ is equal to
$a=a_1\cdot...\cdot a_m$. We denote by ${\cal H}_{(d_1,...,d_m)}$
the space of mappings of the form \RE{defh1}. Each $g\in {\cal
H}_{(d_1,...,d_m)}$ depends on $k$ complex and therefore on $2k$
real variables for some positive integer $k$. Note that $k$ is
determined by $d_1,...,d_m$. For all $d\geq 2$ we define
\begin{displaymath}
{\cal H}_d=\bigcup {\cal H}_{(d_1,...,d_m)},
\end{displaymath}
where the union is taken over all $(d_1,...,d_m)$ with $d_i\geq 2$
such that $d=\nolinebreak
d_1\nolinebreak\cdot\nolinebreak.\nolinebreak.\nolinebreak.\nolinebreak\cdot\nolinebreak
d_m$. Since dynamical properties are invariant under conjugacy,
each polynomial automorphism of $\bC^2$ with nontrivial dynamics
is represented in ${\cal H}_d$ for some positive integer $d$. In
the remainder of this paper we will restrict our attention to
mappings in ${\cal H}_d$.
\\[0.3cm] We denote by $Hyp_d$ the space of hyperbolic mappings in
${\cal H}_d$. A mapping $g\in {\cal H}_d$ is called hyperbolic if
its Julia set $J$ is a hyperbolic set of $g$. That is,  the
tangent bundle of $J$ allows an invariant splitting $E=E^s\oplus
E^u$ such that the tangent map is uniformly contracting on $E^s$
and uniformly expanding on $E^u$ (see \cite{KH} for further
details). The results of \cite{BS1} imply that in the hyperbolic
situation $J$ has index $1$; that is $\dim_\bC E^{u/s}_p =1$ for
all $p\in J$. The nonwandering set of $g$ is equal to the union of
$J$ with finitely many sinks. Furthermore $g$ is an Axiom A
diffeomorphism and $J$ is a basic set of $g$. The most important
feature of hyperbolic sets is that we can associate with each
point $p$ the local unstable/stable manifold
$W^{u/s}_\epsilon(p)$,  varying continuously with $p$. We will
also denote by $W^{u/s}(p)$ the (global) unstable/stable manifold
of $p$. In the case of mappings in $Hyp_d$ the (local)
unstable/stable manifolds are in fact complex manifolds of complex
dimension one. Finally we note that $Hyp_d$ is open.\\[0.3cm] We
will work both with the Hausdorff and the topological dimension of
a set $A$. For the Hausdorff dimension, denoted by $dim_H A$, see
\cite{M}; and for the topological dimension, denoted by $dim_{top}
A$, see \cite{HW}. For any $A$ we have $dim_{top} A\leq dim_H A$.
 \section{Topological dimension of $J$}
In this section we establish some results about the topological
dimension of the Julia set of a polynomial automorphism of
$\bC^2$. In particular we show that in the hyperbolic case the
topological dimension of the Julia set is given by an intersection
formula, see Theorem \ref{thintsec}. We also conclude that if $g$
is  hyperbolic, then the Julia set cannot be locally connected.
\\[0.3cm] We start with a basic Lemma.
\begin{lemma}
Let $g\in Hyp_d$. Then
$t^{u/s}_{top}=dim_{top}W^{u/s}_\epsilon(p)\cap J$ does not depend
on $p\in J$.
\end{lemma}
{\it Proof. } Without loss of generality we only consider the
unstable slice. Let $p,q\in J$ such that $|p-q|$ is small. Then it
is a general result about hyperbolic sets that  the holonomy
mapping maps $W^u_\epsilon (p)\cap J$ homeomorphically to
$W^u_\epsilon (q)\cap J$. Thus $t^u_{top}$ is locally constant. It
is shown in \cite{BS1} that $g|_J$ is topologically transitive. We
conclude that $t^u_{top}$ is independent of $p\in J$. \eproof We
refer to $t^{u/s}_{top}$ as the topological dimension of the
unstable/stable slice. The topological dimension of these slices
is related to the topological dimension of $J^\pm$.
\begin{theorem}\label{thtopdim}
Let $g\in Hyp_d$. Then $dim_{top} J^\pm = t^{u/s}_{top} +2.$
\end{theorem}
{\it Proof. } Without loss of generality we show the result only
for the unstable slice. Let $p\in J$. The Stable Manifold Theorem
implies that there exists a continuous mapping
\begin{displaymath}
\Psi_p:\overline{W}^u_\epsilon(p)\cap J\times
\overline{D}(0,\epsilon)\to \bC^2
\end{displaymath}
such that
$\Psi_p(q,\overline{D}(0,\epsilon))=\overline{W}^s_\epsilon(q)$
and that $\Psi_p(q,.)$ is injective for all $q\in
\overline{W}^u_\epsilon(p)\cap J$ (see \cite{W1} for further
details). Here $\overline{W}^{u/s}_\epsilon(p)$ and
$\overline{D}(0,\epsilon)$ denote the closed local unstable/stable
manifold of $p$ and the closed disk in $\bC$ with center $0$ and
radius $\epsilon$ respectively. Since $J$ is a hyperbolic set for
$g$, the mapping $g|_J$ is expansive. Thus we can assure by making
$\epsilon$ smaller if necessary that $\Psi_p$ is injective and
therefore a homeomorphism onto its image. Therefore
\begin{equation}
dim_{top}\left(\bigcup_{q\in \overline{W}^u_\epsilon(p)\cap J}
\overline{W}^s_\epsilon(q)\right)=t^u_{top} + 2.
\end{equation}
As in the proof of Theorem 4.1 of \cite{W3}, there exist
$p_1,...,p_n\in J$ and $\epsilon_1,...,\epsilon_n>0$ such that for
\begin{equation}
\epsilon=\min\left\{\frac{\epsilon_1}{2},...,\frac{\epsilon_n}{2}\right\}
\end{equation}
and for all $p\in J$  we have that $\overline{W}^s_{\epsilon}(p)$
is contained in $\overline{W}^s_{\epsilon_k}(q)$ for some $q\in
\overline{W}^u_{\epsilon_k}(p_k)\cap J$ and some
$k\in\{1,...,n\}$. Note that this follows from the fact that $J$
is a hyperbolic set with a local product structure. We conclude
that
\begin{equation}
dim_{top} \left(\bigcup_{p\in
J}\overline{W}^s_{\epsilon}(p)\right)=t^u_{top}+2.
\end{equation}
It is a result of Bedford and Smillie \cite{BS1} that
$W^{s}(J)=J^+$. Therefore, we can conclude by  Proposition 3.10 of
\cite{B1} that
\begin{equation}
\bigcup_{p\in J}W^{s}(p)=J^+.
\end{equation}
On the other hand, we have
\begin{equation}
\bigcup_{n\in\bN}g^{-n}\left(\bigcup_{p\in
J}\overline{W}^{s}_{\epsilon}(p)\right)=\bigcup_{p\in J}W^{s}(p).
\end{equation}
Hence
\begin{equation} dim_{top} J^+= 2+t^u_{top}
\end{equation}
and the proof is complete. \eproof The reason for considering
closed local unstable/stable manifolds in the proof of Theorem
\ref{thtopdim} is that some basic properties for the topological
dimension only hold in general  for closed sets (see \cite{HW}).
\begin{corollary}\label{cordimtop01}
Let $g\in Hyp_d$. Then $t^{u/s}_{top}\in \{0,1\}$.
\end{corollary}
{\it Proof. } A subset of $\bR^n$ has topological dimension $n$ if
and only if it has non-empty interior (see \cite{HW}). This
implies $\dim_{top}J^\pm\leq 3$. Therefore the result follows
immediately from Theorem \ref{thtopdim}.\eproof { Remark. } Since
the topological dimension of a set is a lower bound for its
Hausdorff dimension, the result of Corollary \ref{cordimtop01}
also follows from  $dim_H W^{u/s}_\epsilon(p)\cap J< 2$ (see
\cite{W3}).\\[0.3cm] Let $g\in Hyp_d$. By \cite{BS1}, the Julia
set $J$ has a local product structure; this implies
\begin{equation}\label{eqprostr}
dim_{top} J= t^u_{top}+ t^s_{top}.
\end{equation}
Therefore, Corollary \ref{cordimtop01} yields
\begin{equation}\label{eq012}
dim_{top}J\in \{0,1,2\}.
\end{equation}
 Similarly as it was done in \cite{W3} for the Hausdorff dimension, Theorem \ref{thtopdim} and equation \RE{eqprostr}
 imply
an intersection formula for the topological dimension of the Julia
sets.
\begin{theorem}\label{thintsec}
Let $g\in Hyp_d$. Then
\begin{equation}\label{eqtd}
dim_{top} J =dim_{top} J^+ + dim_{top} J^- -4.
\end{equation}
\end{theorem}
We refer to equation \RE{eqtd} as the intersection formula for the
topological dimension.\\ As $J$ is compact, its topological
dimension exceeds $0$ if and only if $J$ is not totally
disconnected. This implies that a mapping $g\in {\cal H}_d$ with a
connected Julia set $J$ provides an example for $dim_{top} J >0$.
In the following result we present another criterion in terms of
the existence of particular Fatou components.
\begin{theorem}\label{thtop1}
Let $g\in {\cal H}_d$ and $|a|<1$. Let $p\in J$ be a saddle point
of $g$. Assume there exists a periodic connected component $C$ of
$int K^+$ such that $C\cap J^-\not=\emptyset$. Then $dim_{top}
W^u_{\epsilon}(p)\cap J \geq 1$.
\end{theorem}
{\it Proof. } Without loss of generality we  assume that $p$ is a
fixed point and $C$ has period $1$. Since  $C\cap
J^-\not=\emptyset$, we may follow by Proposition 7 of \cite{BS2}
that $\partial C = J^+$. We have  $W^u_{\epsilon}(p)\subset J^-$.
This implies $W^u_{\epsilon}(p)\cap J = W^u_{\epsilon}(p)\cap
J^+$. Together we obtain that
\begin{equation}
W^u_{\epsilon}(p)= (W^u_{\epsilon}(p)\cap C) \cup
(W^u_{\epsilon}(p)\cap J)\cup (W^u_{\epsilon}(p)\setminus
\overline{C})
\end{equation}
is a disjoint union. We will show  that none of the sets in this
union is empty. The set $W^u_{\epsilon}(p)\cap J$ is not empty
because it contains $p$. Since $C$ is $g$-invariant and
$\overline{W^u(p)}=J^-$ (see \cite{BS2}), we conclude that
$W^u_\epsilon(p)\cap C\not=\emptyset$. We have
$W^u_\epsilon(p)\cap C, W^u_\epsilon(p)\cap J\subset K$. This
implies that
\begin{equation}
\bigcup_{n\in \bN} g^n(W^u_\epsilon(p) \cap \overline{C})
\end{equation}
is bounded. The unstable manifold $W^u(p)$ is a holomorphic copy
of $\bC$ and is therefore unbounded. Thus
$W^u_\epsilon(p)\setminus \overline{C}\not=\emptyset$. The local
unstable manifold $W^u_\epsilon(p)$ is a topological disk and
therefore cannot be separated  by a subset of topological
dimension $0$  (see \cite{HW}). We conclude that $dim_{top}
W^u_{\epsilon}(p)\cap J \geq 1$.\eproof  { Remarks. }\\ Note that
no hyperbolicity is required in Theorem \ref{thtop1}.  Bedford and
Smillie classified in \cite{BS2}
 recurrent connected components of $int K^+$ in the case $|a|<1$. They
showed that such a component $C$ is periodic and must  either be a
basin of attraction or the stable set of a rotational domain
(Siegel disk, Herman ring), and that in all  these cases  $C\cap
J^-\not=\emptyset$ holds. Therefore Theorem \ref{thtop1} applies.
Note that Theorem \ref{thtop1} also holds in the case of a
non-recurrent periodic connected component $C$ of $int K^+$ which
satisfies $C\cap J^-\not=\emptyset$. If $g\in Hyp_d$, then by
\cite{BS1} the only possibility for a connected component of $int
K^+$   is a basin of attraction.
\\[0.3cm] Finally, we show that  $J$ cannot be locally connected.
\begin{theorem}\label{thloccon}
Let $g\in Hyp_d$. Then the Julia set $J$ is not locally connected.
\end{theorem}
{\it Proof. } Let us assume $J$ is locally connected and $p\in J$
is a saddle fixed point. The local product structure of $J$
implies that there exists a connected neighborhood $U\subset J$ of
$p$ such that $U$ is homeomorphic to $(W^u_\epsilon(p)\cap J)
\times (W^s_\epsilon(p)\cap J)$. Therefore each of the sets,
$W^u_\epsilon(p)\cap J$, $ W^s_\epsilon(p)\cap J$, is connected
and has topological dimension equal to $1$. It is shown in
\cite{BS1} that $g|_J$ is topologically mixing. This implies
\begin{equation}
\overline{\bigcup_{n\in \bZ} g^n(U)} = J.
\end{equation}
We conclude that $J$ is connected, and therefore $g$ is an
unstably connected mapping (see \cite{BS6}). Thus, we may conclude
by Corollary 5.9 of \cite{BS7} that $W^s_\epsilon(p)\cap J$ is a
Cantor set. This is a contradiction to  $dim_{top}
W^s_\epsilon(p)\cap J=1$ which proves the desired result. \eproof
\section{Topological pressure, Bowen- Ruelle formula and fractal Julia
sets} In this section we apply the Bowen- Ruelle formula to obtain
 information about the Hausdorff dimension of the Julia set of a
hyperbolic polynomial automorphism of $\bC^2$. In particular, we
use properties of topological pressure to show that the Hausdorff
dimensions of the unstable and stable slices are related to each
other by the Jacobian determinant of the mapping (see Corollary
\ref{corneu}). This relation is then applied in Theorem
\ref{thfrac} to show that   the Hausdorff dimension of the Julia
set always exceeds its topological dimension. \\[0.3cm] First we
recall the concept of topological pressure and some of its
properties. The main references for the following are \cite{R} and
\cite{W}.

Let $(X,d)$ be a compact metric space and $T:X\to X$ a continuous
mapping. For $n\in\bN$ we define a new metric $d_n$ on $X$ given
by $d_n(x,y)=\max\limits_{i=0,...,n-1} d(T^i(x),T^i(y))$. Let
$\epsilon>0$. A set $F\subset X$ is called
$(n,\epsilon)$-separated with respect to $T$ if
$d_n(x,y)<\epsilon$ implies $x=y$ for $x,y\in F$. For all
$(n,\epsilon)\in \bN\times \bR^+$ let $F_n(\epsilon)$ be a maximal
$(n,\epsilon)$-separated set with respect to $T$ (in the sense of
inclusion). We denote by $C(X,\bR)$ the Banach space of all
continuous functions from $X$ to $\bR$. The topological pressure
of $T$ is a mapping $P(T,.)$ from $C(X,\bR)$ to $\bR\cup
\{\infty\}$ defined by
\begin{equation}\label{defpre}
P(T,\varphi)=\lim_{\epsilon \to
0}\limsup_{n\to\infty}\frac{1}{n}\log\left(\sum_{x\in
F_n(\epsilon)} \exp\left( \sum_{i=0}^{n-1}\varphi \circ
T^i(x)\right)\right).
\end{equation}
If $T$ is expansive, then we can omit the limit $\epsilon\to 0$
and take the value for a fixed $\epsilon>0$ with the property that
$2\epsilon$ is a constant of expansivity  for $T$. Note that
$P(T,0)= h_{top}(T)$, where $h_{top}(T)$ denotes the topological
entropy of $T$. We define $M(X,T)$ to be the space of all
$T$-invariant Borel probability measures on $X$. The variational
principle provides the formula
\begin{equation}\label{eqvarpri}
P(T,\varphi)= \sup_{\mu\in M(X,T)} \left(h_\mu(T)+\int_X \varphi
d\mu\right),
\end{equation}
 where $h_\mu(T)$ denotes the measure theoretic entropy of $T$
with respect to $\mu$. In the sequel we assume that
$h_{top}(T)<\infty$. The topological pressure has the following
properties. \\\beha
\item  The topological pressure  is a convex function.
\item If $\varphi$ is a strictly negative function, then the
function $t\mapsto P(T,t\varphi)$ is strictly decreasing.
 \item For all $\varphi\in C(X,\bR)$ we have
 \begin{equation}\label{eqpab}
 h_{top}(T)+\min \varphi \leq P(T,\varphi)\leq h_{top}(T)+\max
 \varphi.
 \end{equation}
 \behe
 In many hyperbolic
systems topological pressure is related to dimension. We will
discuss this relation in the case of polynomial automorphisms of
$\bC^2$. Let $g\in Hyp_d$. We define
\begin{displaymath}
\phi^{u/s}:J\to \bR\qquad p\mapsto \log ||Dg(p)|_{E^{u/s}_p}||
\end{displaymath}
and  the unstable/stable pressure function
\begin{displaymath}
P^{u/s}:\bR\to\bR\qquad t\mapsto P(g|_J,\mp t\phi^{u/s}).
\end{displaymath}
 The following result due to Verjovsky and Wu \cite{VW} is
crucial for our further presentation.
\begin{theorem}\label{thVW}
Let $g\in Hyp_d$ and $p\in J$. Then $t^{u/s}=dim_H
W^{u/s}_\epsilon(p)\cap J$ does not depend on $p\in J$.
Furthermore $t^{u/s}$ is given by the unique solution of
\begin{equation}\label{eqBR}
P^{u/s}(t)=0.
\end{equation}
\end{theorem}
Equation \RE{eqBR} is called Bowen- Ruelle formula. We refer to
$t^{u/s}$ as the Hausdorff dimension of the unstable/stable slice.
Theorem \ref{thVW} we can be applied  to derive rough bounds for
the Hausdorff dimension of $J$. In comparison to equation
\RE{eqBR} the dynamical  meaning of these bounds is more
transparent.
\begin{theorem}Let $g\in Hyp_d$.
We define
\begin{displaymath}
\begin{array}{l}
\overline{s}=
\lim\limits_{n\to\infty}\frac{1}{n}\log\left(\max\{||Dg^n(p)|_{E^u_p}||:p\in
J\}\right)\\[0.3cm] \underline{s}=
\lim\limits_{n\to\infty}\frac{1}{n}\log\left(\min\{||Dg^n(p)|_{E^u_p}||:
p\in J\}\right).  \end{array} \end{displaymath} Then
\begin{equation}\label{eqpromo}
\ \left(\frac{1}{\overline{s}}+\frac{1}{\overline{s}-\log
|a|}\right)\log d\leq dim_H J \leq
\left(\frac{1}{\underline{s}}+\frac{1}{\underline{s}-\log
|a|}\right)\log d.
\end{equation}
\end{theorem}
{\it Proof. }  Since $E^{u/s}_p$ is of complex dimension one, we
have
\begin{equation}\label{eqmul}
||Dg^n(p)|_{E^{u/s}_p}||= \prod_{i=0}^{n-1}
||Dg(g^i(p))|_{E^{u/s}_{g^i(p)}}||
  \end{equation}
for all $n\in \bN$ and all $p\in J$.  Therefore, the limits
defining $\underline{s}$ and $\overline{s}$ exist (see \cite{W1}
for further details). Since $g$ is hyperbolic, we have
$\underline{s}> 0$. Let us consider the upper bound in inequality
\RE{eqpromo}. For $n\in\bN$ we define
\begin{displaymath}
\phi^u_n=\frac{1}{n} \log ||Dg^n|_{E^u}||. \end{displaymath}  It
follows from equation \RE{eqmul} that
\begin{equation}\label{eqint}
\int \phi^u d\mu =  \int \phi^u_n d\mu
\end{equation}
for all $n\in\bN$ and all $\mu\in M(J,g|_J)$  (see \cite{W3}). We
conclude by equation \RE{eqvarpri} and Theorem \ref{thVW} that
\begin{equation}
P(g|_J, -t^u \phi^u_n)=0
\end{equation}
for all $n\in \bN$. Let us recall that $h_{top}(g|_J)=\log d$ (see
\cite{BS3}). Therefore equation \RE{eqpab} yields
\begin{equation}
\log d- t^u{ \textstyle
\frac{1}{n}}\log\left(\min\{||Dg^n(p)|_{E^u_p}||: p\in
J\}\right)\geq 0.
\end{equation}
Taking the limit implies $t^u\leq \frac{\log d}{\underline{s}}$.\\
In \cite{W1} we have shown that
\begin{equation}
\underline{s}-\log|a|=\lim\limits_{n\to\infty}
\frac{1}{n}\log\left(\min\{||Dg^{-n}(p)|_{E^s_p}||: p\in
J\}\right).
\end{equation}
Therefore the same argumentation as above applied to $g^{-1}$
yields $t^s\leq \frac{\log d}{\underline{s}-\log |a|}$. Applying
\begin{equation}
dim_H J= t^u+t^s \end{equation} completes the proof for the upper
bound in inequality \RE{eqpromo}. Analogously we obtain the proof
for the lower bound in inequality \RE{eqpromo}.\eproof { Remark. }
 Friedland and Ochs   provided in \cite{FO} a
formula for the Hausdorff dimension of the Julia set in terms of
entropies and Lyapunov exponents of invariant probability measures
(see also \cite{W3}). This result can also be applied to derive
lower and upper bounds for the Hausdorff dimension of $J$ similar
to those in equation \RE{eqpromo}.\\[0.3cm]
 Now
we introduce a natural generalization of bidisks in $\bC^2$.
\begin{definition}
Let $E_1,E_2\subset \bf C^2$ be complex lines through the origin
such that $\bC^2= E_1 \bigoplus E_2$. Let $p\in \bC^2$ and
$r=(r_1,r_2)$, where $r_1,r_2>0$. We define the  $(E_1,
E_2)$-bidisk with center $p$ and radius $r$    by
\begin{displaymath}
P_{E_1,E_2}(p,r_1,r_2)=\{p+e_1+e_2:e_1\in E_1,e_2\in E_2,
|e_1|<r_1,|e_2|<r_2\}.
\end{displaymath}
Furthermore $P_{E_1,E_2}(r_1,r_2)=P_{E_1,E_2}(0,r_1,r_2)$.
\end{definition}
We need the following  Lemma.
\begin{lemma}\label{lem10}
Let $g\in Hyp_d$. Then there exist $C_1,C_2>0$ such that for all
$p\in J$ and all $r_1,r_2>0$ we have
\begin{equation}
C_1r_1^2 r_2^2\ \leq vol(P_{E^s_p,E^u_p}(r_1,r_2))\leq C_2 r_1^2
r_2^2.
\end{equation}
\end{lemma}
{\it Proof. } Let $(e_p^s)_{p\in J}$ and $(e_p^u)_{p\in J}$ be
continuous families in $\bC^2$ such that  $ e^s_p \in E^s_p, e^u_p
\in E^u_p$  and $|e^s_p|, |e^u_p|=1$ for all $p\in J$. For $p\in
J$ we define a linear mapping
\begin{displaymath}
A_p:\bC^2\to\bC^2\qquad ( z_1,z_2) \mapsto z_1 e^s_p + z_2 e^u_p.
\end{displaymath}
Since $J$ is a hyperbolic set we have that
\begin{displaymath}
C_1=\pi^2 \min_{p\in J}\{|\det\, A_p|^2\} \qquad C_2= \pi^{2}
\max_{p\in J}\{|\det\, A_p|^{2}\}
\end{displaymath}
are well-defined and $0<C_1\leq C_2<\infty$.\\ Obviously
$P_{\bC\times \{0\},\{0\}\times \bC}(r_1,r_2)$ is mapped
bijectively by $A_p$ to $P_{E^s_p,E^u_p}(r_1,r_2)$. Let $r_1,r_2
>0$ and $p\in J$. Then
\begin{equation}
\begin{array}{lll}
C_1 r_1^2 r_2^2 &=& C_1 \pi^{-2} vol(P_{\bC\times
\{0\},\{0\}\times \bC}(r_1,r_2))\\[0.15cm] &=& C_1 \pi^{-2}|\det\,
A_p|^{-2}vol(P_{E^s_p,E^u_p}(r_1,r_2))\\[0.15cm]
&\leq&vol(P_{E^s_p,E^u_p}(r_1,r_1))\\[0.15cm] &=& |\det\,
A_p^{-1}|^{-2} vol(P_{\bC\times \{0\},\{0\}\times
\bC}(r_1,r_1))\\[0.15cm] &\leq&  C_2 r_1^2 r_2^2.
\end{array}
\end{equation}
This completes the proof. \eproof As an immediate consequence we
obtain the following.
\begin{lemma}\label{lempro}
Let $g\in Hyp_d$. Then there exists $C_1,C_2>0$ such that
\begin{equation}
C_1 |a|^n\leq ||Dg^n(p)|_{E^s_p}||\cdot ||Dg^n(p)|_{E^u_p}||\leq
C_2|a|^n
\end{equation}
for all $p\in J$ and all $n\in \bN$.
\end{lemma}
{\it Proof. } We have
\begin{equation}
Dg^n(p)(P_{E^s_p,E^u_p}(1,1))=
P_{E^s_{g^n(p)},E^u_{g^n(p)}}(||Dg^n(p)|_{E^s_p}||,||Dg^n(p)|_{E^u_p}||).
\end{equation}
Therefore the result follows  from Lemma \ref{lem10}. \eproof Next
we  show that the unstable and stable pressure functions are
related to each other by the Jacobian determinant of the mapping
$g$.
\begin{proposition}\label{thtp}
Let $g\in Hyp_d$   and $t\geq 0$. Then $P^u(t)=P^s(t)-t\log |a|.$
\end{proposition}
{\it Proof. } Let $\epsilon>0$ such that $2\epsilon$ is a constant
of expansivity  for $g|_J$. Let $C_1,C_2$ be the constants in
Lemma \ref{lempro}. For all $n\in\bN$ let $F_n(\epsilon)$ be  a
maximal $(n,\epsilon )$-separated set with respect to $g|_J$. Let
$t\geq 0$. Then
\begin{equation}
\begin{array}{ll}
& \frac{1}{n}\log\left(\sum_{p\in
F_n(\epsilon)}\exp\left(\sum_{i=0}^{n-1}-t\phi^u\circ
g^i(p)\right)\right)\\[0.3cm] =& \frac{1}{n}\log\left(\sum_{p\in
F_n(\epsilon)}||Dg^n(p)|_{E^u_p}||^{-t}\right)\\[0.3cm] \leq &
\frac{1}{n}\log\left(\sum_{p\in F_n(\epsilon)}\left(C_1 |a|^n
||Dg^n(p)|_{E^s_p}||^{-1}\right)^{-t}\right)\\[0.3cm] =&
\frac{1}{n}\log\left(C_1^{-t} |a|^{-tn}\sum_{p\in F_n(\epsilon)}
||Dg^n(p)|_{E^s_p}||^{t}\right)\\[0.3cm] =& \frac{-t}{n}\log C_1
-t\log |a| + \frac{1}{n}\log\left(\sum_{p\in
F_n(\epsilon)}\exp\left(\sum_{i=0}^{n-1}t\phi^s\circ
g^i(p)\right)\right).
\end{array}\end{equation}
Taking the upper limit   implies $P^u(t)\leq P^s(t)-t\log |a|.$
The opposite  inequality follows analogously. \eproof{ Remark. }
Note that this result is essentially based on the facts that the
stable and unstable spaces are of complex dimension one and that
$g$ has constant Jacobian determinant.\\[0.3cm] Proposition
\ref{thtp} implies an inequality which relates the Hausdorff
dimensions of the unstable and stable slice.
\begin{corollary}\label{corneu}
Let $g\in Hyp_d$. Then
\begin{equation}\label{eqj-}
t^s\leq \frac{t^u \log d}{\log d -t^u\log|a|}.
\end{equation}
In particular, $t^s<t^u$ for $|a|<1$, and $t^s=t^u$ for $|a|=1$.
\end{corollary}
{\it Proof.} Let us recall our standing assumption $|a|\leq 1$. It
is well-known that $P^u(0)=P^s(0)=h_{top}(g|_J)=\log d$. On the
other hand, we deduce from Proposition \ref{thtp} that
$P^s(t^u)=t^u\log |a|$. Since $t\mapsto P^s(t)$ is a convex
function, its graph lies below the line segment joining $(0,\log
d)$ and $(t^u,t^u\log |a|)$. Therefore application of Theorem
\ref{thVW} implies equation \RE{eqj-}. The rest is trivial.\eproof
\begin{corollary} Let $g\in Hyp_d$. Then\\[0.2cm]
{\rm i)} If $|a|<1$, then \begin{equation}\label{eqende} dim_H J^-
\leq \frac{t^u \log d}{\log d -t^u\log|a|}+2  < dim_H J^+.
\end{equation}
{\rm ii)}
 If $|a|=1$, then $dim_H
J^+=dim_H J^-$.
\end{corollary}
{\it Proof. } We showed in \cite{W3} that $dim_H J^\pm = t^{u/s}
+2.$ Therefore, the result follows from Corollary \ref{corneu}.
\eproof Inequality \RE{eqj-} can also be applied to provide a
criterion for occurrence of a totally disconnected stable slice.
\begin{corollary}\label{corcan}
Let $C$ be a connected component of $Hyp_d$. Assume there exists
$g\in C$ with $|a|\leq d^{-1/2}$. Then the stable slice of every
$f\in C$ is a Cantor set.
\end{corollary}
{\it Proof. } The stable slice is a Cantor set if and only if its
topological dimension is equal to $0$. By the stability of
hyperbolic sets we conclude that the stable slice has  constant
topological dimension  in $C$. Therefore, it is sufficient to show
the result for $g$. We have $t^u<2$  (see \cite{W3}). Thus
equation \RE{eqj-} implies
\begin{equation}
t^s<\frac{2\log d}{\log d -2\log |a|}.
\end{equation}
From $|a|\leq d^{-1/2}$ we conclude that $t^s<1$.  \eproof {
Remark. } In particular we have shown that if $g\in Hyp_d$ and
$dim_{top} J=2$, then $|a|$ has to be sufficiently close to $1$.
Let us mention that  the stable slice of $g\in Hyp_d$ is also
totally disconnected under the assumption that $J$ is connected
(see \cite{BS7}).
\begin{corollary}\label{cor1+s}
Let $g\in Hyp_d$ and $dim_{top} J \geq 1$. Then
\begin{equation}\label{eq1+s}
dim_H J \geq 1 + t^s>1.
\end{equation}
\end{corollary}
{\it Proof.}  i) $t^u_{top}=1$: Since the topological dimension of
a set is a lower bound of its Hausdorff dimension, we obtain
$t^u\geq 1$. We have $t^s>0$. Therefore equation \RE{eq1+s}
follows from $dim_H J = t^u+t^s$.\\ ii) $t^s_{top}=1$: Corollary
\ref{corneu} implies $1=t^s_{top}\leq t^u$ and therefore equation
\RE{eq1+s} follows as in i). \eproof Note that in view of Theorem
\ref{thtop1}, Corollary \ref{cor1+s} is in particular applicable
if $g\in Hyp_d$ has an attracting periodic point.
\\[0.3cm] Finally, we combine our observations  to obtain the main  result of this section.
\begin{theorem}\label{thfrac}
Let $g\in Hyp_d$. Then $ dim_{top} J < dim_H J.$
\end{theorem}
{\it Proof.} Without loss of generality $|a|\leq 1$.\\ If
$dim_{top}J=0$, the inequality follows from the fact that the
Hausdorff dimension of $J$ is strictly positive (see \cite{VW},
\cite{BS3} and \cite{W1}). If $dim_{top} J=1$, the estimate
follows from Corollary \ref{cor1+s}.\\ It remains to consider the
case $dim_{top} J = 2$.\\ If $|a|<1$, then Corollary
\ref{cordimtop01}, Corollary \ref{corneu}  and equation
\RE{eqprostr} imply
\begin{equation} dim_{H}J=t^u+t^s= 2t^s +(t^u-t^s)\geq 2
t^s_{top} + (t^u-t^s)= 2+(t^u-t^s)>2.
\end{equation}
To complete the proof we have to consider the case $dim_{top}J =
2$ and $|a|=1$. Assume there exists $g\in Hyp_d$ such that $dim_H
J=2=dim_{top} J$ and $|a|=1$. This implies that
$t^u=t^u_{top}=t^s=t^s_{top}=1$. Since $g\in Hyp_d$, there exist
generalized  H\'enon mappings $g_1,..., g_m$ such that
$g=g_1\circ...\circ g_m$. Each $g_i$ has the form
\begin{displaymath}
g_i(z,w)=(w,P_i(w)+a_iz),
\end{displaymath}
where $P_i$ is a complex polynomial of degree at least $2$ and
$a_i$ is a non-zero complex number. Note that $a=\prod_{i=1}^m
a_i$. For a fixed $0<\delta<1$ let us consider
$\left(g_\epsilon\right)_{\epsilon\in (-\delta,\delta)}\subset
{\cal H}_d$ defined by
\begin{displaymath}
g_\epsilon=g_{1,\epsilon}\circ g_2\circ...\circ g_m,
\end{displaymath} where
\begin{displaymath}
g_{1,\epsilon}(z,w)=(w,P_1(w)+(1-\epsilon)a_1 z).\end{displaymath}
Note that $\det\, Dg_\epsilon=(1-\epsilon)a$ and $g_0=g$. Let $C$
denote the connected component of $Hyp_d$ containing $g$. Since
$Hyp_d$ is open we can assure by making $\delta$ smaller if
necessary that $g_\epsilon \in C$ for $|\epsilon|<\delta$. This
implies that $\left(g_\epsilon\right)_{\epsilon\in
(-\delta,\delta)}$ is a real analytic family in $Hyp_d$. Thus it
follows from a result of \cite{VW} that $t^u_\epsilon$ and
$t^s_\epsilon$ are real analytic functions of $\epsilon$. The
topological dimension of the unstable/stable slice is constant in
$C$. This implies
 \begin{equation}\label{eqg1}
 t^u_\epsilon, t^s_\epsilon\geq 1
\end{equation}
for all $|\epsilon|<\delta$. Since $t^u_\epsilon$ is real analytic
we can write
\begin{equation}
t^u_\epsilon= 1+\alpha_1 \epsilon+\alpha_2 \epsilon^2+
O(|\epsilon|^3).
\end{equation}
Therefore we conclude from equation \RE{eqg1} that $\alpha_1=0$.
Making $\delta$ smaller if necessary, there exists a positive
constant $\alpha$ such that
\begin{equation}
t^u_\epsilon\leq 1 + \alpha|\epsilon|^2
\end{equation}
for $|\epsilon|<\delta$. Let us now restrict our attention to the
case $\epsilon\in (0,\delta)$. Thus $|\det\,
Dg_\epsilon|=1-\epsilon<1$, that is, $g_\epsilon$ is a volume
decreasing mapping. Application of Corollary \ref{corneu} and an
elementary calculation yield
\begin{equation}
\begin{array}{lll}
\displaystyle
 t^s_\epsilon&\leq & \displaystyle\frac{t^u_\epsilon \log d}{\log
d -t^u_\epsilon \log(1-\epsilon)}\\[0.4cm] &\leq&
\displaystyle\frac{(1+\alpha \epsilon^2) \log d}{\log d -(1+\alpha
\epsilon^2) \log(1-\epsilon)}\\[0.4cm] &\leq &
\displaystyle\frac{\log d +\alpha\epsilon^2 \log d}{\log d
+\epsilon(1+\alpha\epsilon^2)}\\[0.4cm]&\leq&
\displaystyle\frac{\log d +\alpha\epsilon^2 \log d}{\log d
+\epsilon}.
\end{array}
\end{equation}
Therefore, if $\epsilon$  is sufficiently small, then
$t^s_\epsilon<1$. But this is a contradiction to equation
\RE{eqg1} and the proof is complete.\eproof \section{The general
case}
 In this section we present some results about the
dimension of the Julia sets of a general polynomial automorphism
of $\bC^2$. I.e., $g$ is not assumed to be necessarily hyperbolic.
\\[0.3cm] Let $g\in {\cal H}_d$ and $V\subset \bC^2$ compact such
that $K\subset int V$ and $g^{\pm 1}(J^\pm\cap V)\subset J^\pm\cap
V$. It is shown in \cite{BS1} that a closed bidisk $V$ of
sufficiently large radius satisfies this property. We define
\begin{equation}\label{defspm}
 s^\pm=\lim\limits_{n\to\infty}\frac{1}{n}\log\left(\max\{||Dg^{\pm n}(p)||: p\in
J^\pm\cap V\}\right).
\end{equation} The submultiplicativity of the
operator norm guarantees the existence of the limit defining
$s^\pm$. Since all norms in $\bC^n$ are equivalent, the value of
$s^\pm$ is independent of the norm.
\begin{proposition}\label{propunV}
The value of $s^\pm$ is  independent of the choice of $V$.
\end{proposition}
{\it Proof.} We have $W^s(K)=K^+$ and $W^u(K)=K^-$ (see
\cite{BS1}). Therefore, the proof follows by elementary arguments
(see \cite{W1} for details).\eproof For a bounded set $A\subset
\bR^n$ we denote $\overline{dim}_B A$ to be the upper
box-dimension of $A$ (see \cite{M} for the definition). We have
the inequality $dim_H A \leq \overline{dim}_B A$ and equality for
sufficiently regular sets $A$. We will now consider the case when
$g$ is volume decreasing. In this situation it was observed in
\cite{FM} that $K^-$ has Lebesgue measure zero and therefore $K^-
=J^-$. The following theorem provides an even stronger result.
\begin{theorem}\label{th4-}
 Let $g\in {\cal
H}_d$ and $|a|<1$. Then
\begin{equation}\label{eq4-}
\overline{dim}_B K^-\cap V \leq 4-
\frac{2\log\left(|a|^{-1}\right)}{s^-}<4.
\end{equation}
\end{theorem}
{\it Proof.} Note that the real Jacobian determinant of $g^{-1}$
as a mapping of $\bR^4=\bC^2$ is equal to $|a|^{-2}$. Therefore,
the result follows immediately from Theorem 1.1 of
\cite{W2}.\eproof  { Remarks. }\\ i) Since $W^u(K)=K^-$ we can
define an exhaustion $V_k=g^k(V\cap K^-)$ of $K^-$. This implies
that the upper bound in inequality \RE{eq4-} is also an upper
bound for the Hausdorff dimension of $K^-$.\\ ii) In particular
Theorem \ref{th4-} implies that, if $g\in {\cal H}_d$ is not
volume preserving, then the upper box-dimension of $K$ is strictly
smaller than $4$.
\\[0.3cm] We can apply Theorem \ref{th4-} to show that if the
modulus of the Jacobian determinant is small, then the Hausdorff
dimension of $K^-$ is close to $2$. Note that for a hyperbolic
mapping this result already follows from equation \RE{eqende}. Let
$(P_{c_1})_{{c_1}\in {C_1}},...,(P_{c_m})_{{c_m}\in {C_m}}$ be
families of complex polynomials of fixed degree $d_i\geq 2$ such
that $C=C_1\times...\times C_m$ is a compact subset of $\bC^k$ for
some $k\in \bN$. For $a=(a_1,...,a_m)\in (D(0,1)\setminus
\{0\})^m$ and for $c\in C$ we set $g_{a,c}=g_{a_1,c_1}\circ ...
\circ g_{a_m,c_m}$, where $g_{a_i,c_i}(z,w)=(w,P_{c_i}(w)+a_i z)$.
In the following we use the notation $|a|=|a_1\cdot ... \cdot
a_m|$. This implies that $|\det\, Dg_{a,c}|=|a|$ for all $(a,c)\in
(D(0,1)\setminus \{0\})^m\times C$.
\begin{proposition}\label{profs}
For all $\epsilon>0$ there exists $a_0>0$ such that for all
$0<|a|<a_0$ and all $c\in C$ we have $dim_H K^-_{a,c}<2+\epsilon$.
\end{proposition}
{\it Proof.} The compactness of $C$ implies that there exists
$V\subset \bC^2$ compact such that for all $(a,c)\in
(D(0,1)\setminus \{0\})^m\times C$ we have $K_{a,c}\subset int V$
and $g_{a,c}^{-1}(K^-_{a,c}\cap V)\subset K^-_{a,c}\cap V$.
Therefore an elementary calculation implies that there exists
$\alpha >0$ such that
\begin{equation}
\max\left\{||Dg^{-1}_{a,c}(p)||: p\in K^-_{a,c}\cap V\right\}\leq
\alpha |a|^{-1}.
\end{equation}
for all $(a,c)\in (D(0,1)\setminus \{0\})^m\times C$. The result
follows now immediately from Theorem \ref{th4-}. \eproof \\ {
Remark. } Forn{\ae}ss and Sibony \cite{FS} considered complex
H\'enon mappings which are in a particular sense close to a
hyperbolic quadratic polynomial in $\bC$. They showed  the
hyperbolicity of these mappings and that the Hausdorff dimension
of $J^-$ is close to $2$. Proposition \ref{profs} can be
considered as a generalization of this result even to the case of
nonhyperbolic mappings.\\

Finally we apply the H\"{o}lder continuity of the Green function
to obtain a lower bound for the Hausdorff dimension of $J^\pm$.
Let $g\in {\cal H}_d$. Let us recall the following definition
introduced by Hubbard:
\begin{equation}
G^{\pm}(p)= \lim_{n\to\infty}\frac{1}{d^n}\log^+|g^{\pm n}(p)|.
\end{equation}
It is shown in \cite{BS1} that $G^\pm$ is  continuous and
plurisubharmonic on $\bC^2$, and pluriharmonic on $\bC^2\setminus
J^\pm$. The function $G^\pm$ is the Green function of the sets
$K^\pm$ and $J^\pm$, and is  called  Green function of $g^{\pm 1}
$.\\ Fornaess and Sibony showed in \cite{FS} that the Green
function $G^\pm$ is in fact H\"{o}lder continuous and that the
corresponding H\"{o}lder exponent plus 2 provides a lower bound
for the Hausdorff dimension of $J^\pm$. In \cite{W1} we improved
the estimate for the H\"{o}lder exponent and derived the following
result.
\begin{theorem}
Let $g\in {\cal H}_d$ and $s^\pm$ as in equation
\hspace{0.2cm}\mbox{\RE{defspm}}. Then for all $s^\pm_0
> s^\pm$ the Green function $G^\pm$ is H\"{o}lder continuous with H\"{o}lder exponent $\frac{\log
d}{s^\pm_0}$ on every compact subset of $\bC^2$. Furthermore
\begin{equation} dim_H J^\pm\geq 2+\frac{\log d}{s^\pm}.
\end{equation}
\end{theorem}
\vspace{0.6cm}
 {\it Acknowledgements.} I would like to thank Eric
Bedford for his valuable comments on the first draft of the paper.
I also want to thank the referee for carefully reading the
manuscript and for suggesting various improvements.

\end{document}